\documentclass[reqno]{amsart}
\usepackage{amsmath, amsthm, amsfonts, amsopn, amscd, verbatim}

    \theoremstyle{plain}
    \newtheorem{Thm}{Theorem}[section]
    \newtheorem{Prop}[Thm]{Proposition}
    \newtheorem{Lemma}[Thm]{Lemma}
    \newtheorem*{Lemma*}{Lemma}

    \theoremstyle{definition}
    \newtheorem{Def}[Thm]{Definition}
    \newtheorem*{Def*}{Definition}
    \newtheorem{Example}[Thm]{Example}

    \theoremstyle{remark}
    \newtheorem{Remark}{Remark}

    \numberwithin{equation}{section}

    \newcommand{\ie}{\textit{i.e.}}

    \newcommand{\field}[1]{\mathbb{#1}}
    \newcommand{\N}{\field{N}}
    \newcommand{\Z}{\field{Z}}

    \newcommand{\C}{\field{C}}
    \DeclareMathSymbol{\fieldk}{\mathalpha}{AMSb}{"7C} % mathbb k
    \newcommand{\abs}[1]{\lvert#1\rvert}
    
    \newcommand{\norm}[1]{\lVert#1\rVert}
    \newcommand{\inner}[2]{\left\langle#1,#2\right\rangle}
    \DeclareMathOperator{\Tr}{Tr}
    \DeclareMathOperator{\Det}{Det}
    \newcommand{\tensor}{\otimes}
    \DeclareMathSymbol{\normal}{\mathord}{AMSa}{"43}
    \DeclareMathOperator{\Log}{Log}
    \DeclareMathOperator{\Exp}{Exp}
    \DeclareMathOperator{\Id}{Id}
    \DeclareMathOperator{\Image}{Image}
    \DeclareMathOperator{\Ker}{Ker}

    \newcommand{\bs}{\backslash}
    \newcommand{\half}{\frac{1}{2}}
    \newcommand{\g}{\gamma}
    \newcommand{\G}{\Gamma}
    
    \newcommand{\ep}{\varepsilon}

    \newcommand{\Laplace}{\Delta}

    \newcommand{\define}[1]{\emph{#1}}      % Use for defs

    \newcommand{\cover}[1]{\widetilde{#1}}

    \newcommand{\ltwo}{l^{2}}
    \newcommand{\Ltwo}{L^{2}}
    \newcommand{\ctensor}{\hat{\otimes}}
    \newcommand{\VN}[1]{\mathcal{#1}}

    \DeclareMathOperator{\Aut}{Aut}
    \DeclareMathOperator{\vol}{vol}
    \newcommand{\doz}{\partial_1}
    \newcommand{\dz}{\partial_0}
    \newcommand{\Prim}{\mathcal{P}}   % P for classes of ends
    \newcommand{\Ends}{\mathcal{E}}   % E for ends
    \newcommand{\len}{\ell}   % l for length

    \newcommand{\Gx}{\G\!_{\cover{x}}}
    \newcommand{\Ge}{\G\!_{\cover{e}}}
    \newcommand{\Gep}{\G\!_{\ep}}
    \newcommand{\Gepo}{\G_{\ep}^{0}}
    \newcommand{\axis}[1]{X_{#1}}
    \newcommand{\axep}{\axis{\ep}^{+}}
    \newcommand{\lenep}{\len(\ep)}
    \newcommand{\sigep}{\sigma_{\ep}}
    \newcommand{\A}{\VN{A}}
    \newcommand{\Endsn}[1]{\Ends_{n}(#1)}
    \newcommand{\HrhoX}{H_{\rho}(X)}
    \newcommand{\HrhoE}{H_{\rho}(EX)}
    \newcommand{\HrhoV}{H_{\rho}(VX)}
    \newcommand{\HrhoNE}{H_{\rho}(\Endsn{E})}
    \newcommand{\normrho}[1]{\norm{#1}_{\rho}}
    \newcommand{\innerrho}[2]{\inner{#1}{#2}_{\rho}}

\input{boxedeps.tex}
\SetOzTeXEPSFSpecial
\HideDisplacementBoxes

\newcommand{\picta}[3]{
  \begin{figure}[h]   % h=here. Other options: t=top, b=bottom, p=page
    \begin{center}
      \leavevmode
        \ForceWidth{#2}
        \BoxedEPSF{#1.eps}  % other commands: \cBoxedEPSF, \tBoxedePSF
      \caption{#3}
    \end{center}
  \end{figure}
}

\begin{document}

%
% Title Page
%
\title{Zeta Functions Of Discrete Groups Acting On Trees}

\author{Bryan Clair \and Shahriar Mokhtari-Sharghi}
\address{Department of Mathematics,
        Long Island University,
        Brooklyn Campus, 1 University Plaza,
        Brooklyn, NY 11201}
\email{mokhtari@liu.edu}
\address{CUNY Graduate School,
        PhD Program In Mathematics,
        New York, NY 10016}
\email{bclair@gc.cuny.edu}
\date{August 12, 1999}

%\keywords{Cheeger inequality, expanders}
%\subjclass{05c35}

\begin{abstract}
This paper generalizes Bass' work on zeta functions for uniform tree
lattices.
Using the theory of von Neumann algebras, machinery is developed to
define the zeta function of a discrete  
group of automorphisms of a bounded degree tree.
The main theorems relate the zeta
function to determinants of operators defined on edges or
vertices of the tree.

A zeta function associated to a non-uniform tree lattice
with appropriate Hilbert representation is defined. 
Zeta functions are defined for infinite graphs with a
cocompact or finite covolume group action.  
\end{abstract}
\maketitle

%
%  Sections
%

%
% SECTION:
% intro
%
%
%
\section*{Introduction}
The zeta function associated to a finite graph
originated in
work of Ihara~\cite{Ih1,Ih2}, which proves a structure theorem
for torsion-free discrete cocompact subgroups of
$PGL(2,\fieldk_p)$ (where $\fieldk_p$ is a
$p$-adic number field or a field of power series over a finite
constant field). If $\G < PGL(2,\fieldk_p)$ is such a group,
Ihara shows
that $\G$ is in fact a free group and defines a zeta function
associated to $\G$.
An element
$1\neq\gamma\in\G$ is called \define{primitive} if it generates its
centralizer in $\G.$ Let $\Prim(\G)$ denote the set of conjugacy classes
of primitive elements of $\G.$
If $\lambda_1$ and $\lambda_2$ are the eigenvalues of
a representative of $\gamma \in PGL(2,\fieldk_p),$ define
$l(\gamma)= \abs{v_p(\lambda_1\lambda_2^{-1})}$, where $v_p$ is
the normalized
valuation for $\fieldk_p.$
Ihara defines the zeta function attached to $\G$ as:
\[
   Z_\G(u)=\prod_{\gamma\in\Prim(\G)}(1-u^{l(\gamma)} )^{-1}.
\]

In \cite{Ih1}, Ihara extends this definition
to depend on a finite dimensional representation of $\G$ over a field of
characteristic $0$, although in fact the zeta function so
obtained depends only on the character of the given representation.
Ihara's earlier definition is identified
with the zeta
function for the trivial representation of $\G$.
Ihara's zeta function is similar to the Selberg
$\zeta$-function for compact Riemann surfaces.

Ihara's zeta function is defined as an infinite product,
but in fact Ihara proves that it is a rational function.
Let $Z_\G(u,\chi)$ be the Ihara zeta function
with $\chi$ the character of the
representation of $\G$.  Then Ihara's rationality formula states:
\begin{equation}\label{Ihara_ratio}
  Z_\G(u,\chi)=(1-u^2)^s \det(I-Au+pu^2)
\end{equation}
where $s$ is an integer and $A$ is a finite matrix coming out of Ihara's
structure theorem for $\G.$  Using the rationality formula,
one can count the number of
conjugacy classes of primitive elements of $\G.$
In addition, the Ihara
zeta function carries information about the spectral
decomposition of the representation of $G=PGL(2,\fieldk_p)$ on $l^2(\G\bs G).$

Serre \cite[Introduction]{serre} remarked that Ihara's zeta function
can be interpreted as a zeta function of certain
$(p+1)$-regular finite graphs.
In \cite{Ha2}, the authors generalize
some of Ihara's work in the context of finite graphs.
Then, in \cite{Ha1,Ha3}, Hashimoto defines the zeta function of a
finite connected graph $A$ as follows:
For $\rho:\pi_1(A)\to GL(V_\rho)$
any finite dimensional representation of $\pi_1(A),$
the zeta function of $A$ is defined as
\[
  Z(\rho,u) = \prod_{\gamma\in P} \det(I-\rho(\gamma)u^{\len(\gamma)})
  \qquad (\in \C[[u]]),
\]
where $P$ is the set of free homotopy classes of primitive closed paths
in $A$, and $\len(\gamma)$ is the shortest edge length for paths in
the homotopy class of $\gamma$.
To prove a
rationality formula, Hashimoto
assumes that the graph $A$ is bipartite. In
the cases that $A$ is regular or bipartite biregular  
he gives an interpretation for the number
$s$ in \eqref{Ihara_ratio}
in terms of the Euler characteristic of the graph $A.$
The important new step is to express the zeta
function in terms of a rational formula involving edges of the graph.
Specifically,
Hashimoto defines two operators $T_{1,\rho}$ and $T_{2,\rho}$ on the
vector space generated by edges of $A$ and then
proves that the zeta
function of a $A$ is the determinant of an expression involving
$T_{1,\rho}$ and $T_{2,\rho}$.\par

The zeta function of a uniform tree lattice is defined by
H. Bass~\cite{Ba_zeta}.
A \define{uniform tree lattice} consists of a 
tree $X$ and a countable group $\G$ which acts discretely
on $X$ with finite quotient.
Bass' work completes the program of defining zeta functions for
general finite graphs, and furthermore allows for non-free actions
of $\G$ where the quotient is then a graph of finite groups.
To prove the rationality formula, Bass
introduces an operator $T$ on oriented edges of $X$ and
expresses the zeta function in terms of $T$, which greatly simplifies
the situation.

For a general survey and an elementary introduction to
Ihara-type zeta functions, see~\cite{starkter}.

In this paper, we define zeta functions for general
tree lattices bearing appropriate representations which  are  allowed to
be infinite dimensional representations on Hilbert spaces.
This generalization is motivated by an attempt to generalize the work
of Bass to define the zeta function associated to a non-uniform tree
lattice. This generalization
gives rise to a theory of zeta functions for infinite
graphs with symmetry.  Specifically, for
an infinite graph equipped with an action of a discrete group such
that the quotient graph has finite volume. In particular 
one case of interest is that of
non-uniform lattices, where $\G$ and $X$ are as above except that
$\G\bs X$ is infinite but satisfies the finite volume condition:
\[
  \sum_{x \in \G\bs X}\frac{1}{\abs{\G_{\tilde x}}} < \infty,
\]
where $\cover x$ is any representative for $x.$
Before stating precisely the main results, some background is required.

To deal with infinite graphs, we use the machinery of von Neumann
algebras.  The model for
our techniques originated with the $\Ltwo$-Betti numbers, defined by
Atiyah~\cite{atiyah:l2} in 1976.  Since then,
the use of von Neumann algebras to extend classical topological
invariants to non-compact spaces has become fairly well established.
For example, see \cite{luck:l2} for a recent survey.

A von Neumann algebra $\A$ is a *-algebra of bounded operators
on Hilbert space which is closed in the weak (or strong) topology.
The von Neumann algebra $\A$  is equipped with a trace function
$\Tr_{\A} : \A \to \C$. 
A Hilbert $\A$-module $M$ is
a Hilbert space with an $\A$-module structure and an
additional projectivity condition.  The trace extends to a trace for
bounded operators on $M$ which commute with the $\A$ action, and
this trace gives rise to a real-valued dimension theory for Hilbert
$\A$-modules.

Given a formal power series with coefficients which
are trace class operators on a Hilbert $\A$-module,
one can use $\Exp \circ \Tr_{\A} \circ \Log$ to define
a determinant $\Det_{\A}$ as a formal power series over $\C$.
This technique is used in Bass~\cite{Ba_zeta}.
However, the application to von Neumann algebras is apparently
new to this paper and yields a generalization of the classical
real-valued determinant of Fuglede and Kadison~\cite{fk:det}.
We will use the determinant $\Det_{\A}$ to define a zeta
function analogous to the
zeta function for finite graphs.

Let $\G$ be a lattice of a tree $X$ which has uniformly bounded valence.
Let $M$ be a Hilbert $\A$-module, and $\rho$ a representation
from $\G$ to the unitary operators on $M$ which commute with $\A$.

Let $\ep\in\Ends(\G)$ be a hyperbolic end of $X$, that is, $\ep$ is
one of the two ends  of 
 the axis of some hyperbolic element of $\G$.
For $g \in \Gep$(the stabilizer of $\ep$) 
we can define
the (signed) translation distance towards $\ep$, $\tau_{\ep}(g) \in \Z$.
The kernel of $\tau_{\ep}$ is $\Gepo < \Gep$, and define
$\lenep>0$ by the equation $\Image(\tau_\ep)= \lenep\Z$.
Put
\[
  \sigep = \frac{1}{\abs{\Gepo}}
           \sum_{\substack{g\in\Gep \\ \tau_\ep(g)=\lenep}} g
\]

\begin{Def*}
For $\G$ an $X$-lattice and $\rho$ a unitary Hilbert $\A$-representation
of $\G$, define the zeta function of $\G$ by
\begin{equation}
  Z_\rho(u) = \prod_{\ep \in \G\bs\Ends(\G)}
              \Det_{\A}(I-\rho(\sigep)u^{\lenep})\quad \in \C[[u]]
\end{equation}
\end{Def*}
This zeta function is formally quite similar to that of
Bass~\cite{Ba_zeta}.  One crucial difference is that for a
non-uniform lattice there may be infinitely many classes
$\ep \in \G\bs\Ends(\G)$ with the same length, and so it is not
even clear that the coefficients of $Z_{\rho}(u)$ converge.
In fact, conditions on $\rho$ are needed, and will be stated precisely
in Theorem~\ref{thm:mainZT} to follow.

Let $\HrhoE$ denote the
space of $\G$-equivariant functions $f$
from the edges of $X$ to $M$,
which satisfy
\[
  \normrho{f}^{2} =
  \sum_{e \in \G\bs EX}\norm{f(\cover{e})}_{M}^{2} < \infty.
\]
Define the operator $T_{\rho} : \HrhoE \to \HrhoE$ by
\begin{equation*}
    T_\rho f(e) = \sum_{\text{$(e,e_1)$ reduced}} f(e_1).
\end{equation*}

The first theorem of this paper relates the zeta function to the
determinant of $T_{\rho}$.
\begin{Thm}[Bass-Hashimoto Formula]\label{thm:mainZT}
Suppose $X$, $\G$, and $\rho$ satisfy the two finiteness conditions:
\begin{enumerate}
\item $\sum_{e \in \G\bs EX}\dim_{\A}M^{\Ge} < \infty$.
\item For all $n$,
  \begin{equation*}
    \sum_{\substack{\ep \in \G\bs\Ends(\G);\\ \lenep = n}}
        \dim_{\A}M^{\Gepo} < \infty.
  \end{equation*}
\end{enumerate}
Then $Z_{\rho}(u)$ is well defined in $\C[[u]]$, and
\begin{equation}
   Z_{\rho}(u) = \Det_{\A}(I-T_{\rho}u).
\end{equation}
In addition, the power series
$Z_{\rho}(u)$ converges for $u \in \C$, $\abs{u} < \norm{T_{\rho}}$.
\end{Thm}

The two finiteness conditions can intuitively  be interpreted   as
saying  that the 
total volume of the edges of $\G\bs X$ is finite, and the total volume
of loops of length $n$ in $\G\bs X$ is finite for each $n$.  The
subtlety is that even if the stabilizers $\Ge$, $\Gepo$ are large, the
invariant subspaces $M^{\Ge}$ and $M^{\Gepo}$ may not be small.

The next result relates the zeta function to a Laplace operator on
the vertices of $X$.  We define $\HrhoV$ analogous to $\HrhoE$ but
using $\G$-invariant functions on the vertices of $X$.
Put
\begin{equation*}\begin{split}
  \chi_\rho(X) &= \dim_{\A}\HrhoV - \half\dim_{\A}\HrhoE \\
               &= \sum_{x\in VB}\dim_{\A}M^{\Gx} -
                  \half\sum_{e\in EB}\dim_{\A}M^{\Ge}.
\end{split}\end{equation*}
Define the operators $\delta_\rho, Q_\rho : \HrhoV \to \HrhoV$ by
\begin{equation*}
\begin{align}
    \delta_\rho f(x) &= \sum_{\dz e=x}f(\doz e)\\
    Q_\rho f(x)      &= (\deg(x)-1)f(x).
\end{align}
\end{equation*}
Here $\delta_{\rho}$ is the usual adjacency operator.

\begin{Thm}[Bass-Hashimoto-Ihara Formula]\label{thm:mainLap}
When $\sum_{e\in EB}\dim_{\A}M^{\Ge} < \infty$,
\begin{equation*}
  \Det_{\A}(I-T_\rho u) =
  (1-u^2)^{-\chi_\rho(X)}\Det_{\A}(I -  \delta_\rho u +  Q_\rho u^2)
\end{equation*}
where exponentiation is  defined formally using the principal branch
of the logarithm.
\end{Thm}

The remainder of the paper consists of interpreting
Theorem~\ref{thm:mainZT} and Theorem~\ref{thm:mainLap} for important
special cases.

Let $1\to\Lambda\to\G\to\pi\to 1$ be an exact sequence.
The most interesting examples occur when $\A$ is the von Neumann
algebra of the quotient group $\pi = \G/\Lambda$, $\Lambda$ acts
freely on $X$, and $\rho$ is the coset representation of $\G$ on
$M=\ltwo(\pi)$.
In this situation, we may interpret  $Z_{\rho}(u)$ to be
the $\pi$-zeta function
of the infinite graph $Y = \Lambda\bs X$:

\begin{Thm}\label{thm:mainInf}
Suppose $Y$ is a locally finite graph and the discrete group
$\pi$ acts on $Y$ with quotient $B$ having finite volume.
Let $\Laplace_{u} = I - \delta u + Q u^{2}$ acting on $\ltwo(VY)$.
Then
\begin{equation*}
   Z_{\pi}(Y,u) = \prod_{\gamma \in \pi\bs P}
                   \left(1-u^{\len(\gamma)}\right)
                   ^{\frac{1}{\abs{\pi_{\gamma}}}}
                = (1-u^{2})^{-\chi^{(2)}(Y)}\Det_{\pi}(\Laplace_{u}).
\end{equation*}
Here $P$ is the set of
free homotopy classes of primitive closed paths in $Y$.
For $\gamma\in P$, $\len(\gamma)$ is  the length of the shortest
representative of $\gamma$.
The group $\pi_{\gamma}$ is the stabilizer of $\gamma$ under the
action of $\pi$. 
 The Euler characteristic $\chi_{\rho}(X)$ is the $\Ltwo$
Euler characteristic of $\chi^{(2)}(Y)$ as in~\cite{cg:l2coho}, and
can be computed as
\[
  \chi^{(2)}(Y) = \sum_{x \in VB}\frac{1}{\abs{\pi_{x}}} -
                  \half\sum_{e \in EB}\frac{1}{\abs{\pi_{e}}}.
\]
\end{Thm}

The structure of this paper is as follows.  Section one is a review of
graphs and group actions on graphs.  Section two is a treatment of von
Neumann algebras and Hilbert modules, the main novelty of which is
our definition of a complex valued von Neumann determinant.

In the third section, we introduce
the zeta function and prove Theorem~\ref{thm:mainZT}
and Theorem~\ref{thm:mainLap}.
Section four specializes to the main example of a coset
representation and proves Theorem~\ref{thm:mainInf}.
Finally, Section five presents explicit examples.

%
% SECTION:
% tree
%
%
%
\section{Groups and Trees}
\label{group_tree}
This section is a review of graphs and group actions on graphs,
and mainly consists of definitions and notation to be used later.
\begin{Def}
  A \define{graph} $A=(VA,EA)$ consists of a set of
  vertices $VA,$ a set of edges $EA,$ maps
  $\dz,\,\doz :EA\to VA,$ and a fixed point free involution
  $\,\bar{}:EA\to EA$ for which
  $\partial_i \bar{e}=\partial_{1-i}e$ for all $e\in EA.$
\end{Def}
For a graph $A$ and $a\in VA,$ let
$E_0(a)=\dz^{-1}(a)=\{e \in EA \mid \dz e=a\}.$
$A$ is called \define{$k$-regular} if $\abs{E_0(a)}=k$ for all
$a\in VA,$ where
$\abs{E_0(a)}$ is the \define{valence} (degree) of the vertex $a.$
A graph is \define{regular} if it is $k$-regular for some $k.$

The sequence $P=(e_1,e_2,\dotsc,e_n)$ is
a \define{path} of length $n$ from $\dz e_1$ to $\doz e_n$ in $A$ if
$e_i\in EA$
and  $\doz e_i=\dz e_{i+1},$
$1\le i\le n-1.$ For
$a_0\in A$, the empty path from $a_0$ to $a_0$ is denoted $(a_0)$ and
has length $0.$  A path $P$ is \define{closed} if $\dz e_1=\doz e_n$
or if $P$ has length $0.$
$P$ is \define{reduced} if
$e_{i+1}\neq \bar {e_i},\,1\le i\le n-1$ or if $P$ has length $0.$

A nonempty graph $A$ is \define{connected} if any two vertices of
$A$ can
be joined by a path. A graph is a \define{tree} if it is connected and
every reduced path with positive length is not closed. A graph is
\define{locally finite} if the degree of each vertex is finite.
A tree has \define{bounded degree} if there is a uniform bound $d$
on the degrees of vertices.
%In this monograph we always assume $X$ is a bounded degree tree.\par

Let $A$ and $B$ be graphs. A \define{morphism}
$f:A\to B$ is a map taking
$VA\to VB,$ $EA\to EB$ which commutes with the involution
and with $\partial_i.$
A morphism is an \define{automorphism} if it also induces a bijection
between the vertices and edges of $A$ and $B.$ An automorphism of a graph
is called an \define{inversion} if it takes some edge $e$ to $\bar e.$

Let $G\le \Aut(A)$ where $A$ is a connected graph and $\Aut(A)$ is the
group of automorphisms of $A.$ Assume $G$ does not contain any
inversion. The quotient set $G\bs A$  naturally has the structure of a
graph.\par

Now let $X$ be a locally finite tree and $G=\Aut(X)$ the group of
automorphisms of $X.$ Then $G$ is a locally compact group with the compact
open topology coming from the action of $G$ on $X$ ($X$ has the discrete
topology).  The stabilizers of vertices of $X$ are compact and
open. The group $\G \le G$ is discrete if and only if  $\G_x$ is
finite for some
(hence for all) vertices $x\in X.$

Assume that $\G$ does not contain any inversions.
As in~\cite{bass_lattice},
\begin{Def}\label{def:lattice}
The \define{volume} of $\G\bs X$ is defined to be
\[
  \vol(\G\bs X) =\sum_{x \in \G\bs X}\frac{1}{\abs{\G_{\tilde x}}}
\]
where $\tilde x$ is any representative of $x$.
When $\vol(\G\bs X)<\infty$, say that $\G$
has \define{finite covolume} and that
$\G$ is an \define{$X$-lattice}.
In particular $\G$ is called \define{uniform} if $\G\bs X$ is finite and
\define{non-uniform} otherwise.
One can similarly define an \define{edge-volume} by
\[
  \vol(\G\bs EX) =\sum_{e \in \G\bs EX}\frac{1}{\abs{\G_{\tilde e}}}
\]
It is not hard to show that for a bounded degree tree,
   $\vol(\G\bs X) < \infty$ if and 
only if $\vol(\G\bs EX) < \infty$. 
\end{Def}
If $\G$ is an $X$-lattice it is also a lattice of $\Aut(X)$ which
implies that $\Aut(X)$ is unimodular, \ie\ the left Haar measure is
also a right Haar measure.

%
% SECTION:
% von
%
%
\section{Von Neumann Algebras}
The material in this section is mostly classical, although some of
Section~\ref{sec:dets} is new.  References for this material
are \cite{dix:vn} and \cite{lur:kthe}.
\subsection{Traces on von Neumann algebras}
Let $H$ be a separable complex Hilbert space, and let
$B(H,H)$ denote the $C^{*}$-algebra of bounded linear operators on
$H$.  In the \define{weak topology} on $B(H,H)$,
$T_{n} \to T$ if and only if $\inner{T_{n}x}{y} \to \inner{Tx}{y}$
for all $x,y \in H$.  In the \define{strong topology} on $B(H,H)$,
$T_{n} \to T$ if and only if $\norm{T_{n}x} \to \norm{Tx}$ for all
$x \in H$.

\begin{Def}
A \define{von Neumann algebra} is a subalgebra $\A \leq B(H,H)$
which is closed under adjoints and such that $\A$ is closed in the
weak operator topology.
\end{Def}

For a subset $M \subset B(H,H)$, the \define{commutant} of
$M$ is $M^\prime = \{T \in B(H,H) \mid ST = TS, \forall S \in M\}.$
Clearly $M \subset M^{\prime\prime}.$

The double commutant theorem, due to von Neumann, states that the
conditions
\begin{enumerate}
\item $\A^{\prime\prime} = \A$
\item $\A$ is weakly closed
\item $\A$ is strongly closed
\end{enumerate}
are all equivalent.

An element $T \in \A$ is \define{positive} if $T$ is self-adjoint
and $\inner{Tf}{f} \geq 0$ for all $f \in H$.  Let $\A^{+}$ denote
the cone of positive elements.  A map $\Tr : \A^{+} \to [0,\infty]$
is called a \define{trace} if $\Tr(S + T) = \Tr(S) + \Tr(T)$,
$\Tr(\lambda T) = \lambda \Tr(T)$, and
$\Tr(TT^{*}) = \Tr(T^{*}T)$ for all $S,T \in \A$,
$\lambda \in [0,\infty)$.  A trace is \define{finite} if it always
assumes finite values.  A trace is \define{faithful} if $\Tr(T) = 0$
implies $T = 0$ for $T \in \A^{+}$.  A trace is \define{normal}
if for any monotone increasing
net $\{T_{i}\}_{i \in I}$ ($T_{i} - T_{j}$ is positive for $i > j$)
with
$T = \sup T_{i}$, one has $\Tr T = \sup \Tr T_{i}$.

A von Neumann algebra is called \define{finite} if it has a finite,
faithful, normal trace.  A finite trace extends uniquely to a
$\C$-linear map $\Tr_{\A}: \A \to \C$, which satisfies
$\Tr_{\A} ST = \Tr_{\A} TS$ for all $S,T \in \A$.  From now on,
all von Neumann algebras will be finite and equipped with
corresponding trace.

The trace on $\A$ gives an inner product on $\A$ defined by
$\inner{S}{T} =  \Tr(T^{*}S)$.  Let $\hat\A$ denote the Hilbert
completion of $\A$.  $\A$ acts on the left on $\hat\A$,
and we let $B^{\A}(\hat\A)$ denote
the space of bounded, $\A$-equivariant operators on $\hat\A$.
The opposite algebra $\A^{op}$ acts on the right on $\hat\A$,
and the right representation
$R : \A^{op} \to B^{\A}(\hat\A)$ is an isomorphism
(see \cite{dix:vn} for a proof of this fundamental result).
      %more specific refrence is better. Sh.
      % Alas - I can't get ahold of the book.  I agree, though.
Therefore the trace on $\A$ gives rise to a trace on
$B^{\A}(\hat\A)$, also denoted by $\Tr_{\A}$.

\begin{Example}\label{ex:triv}
One can take $\A = \C \subset B(H,H)$, and set $\Tr_{\A}c = c$.
\end{Example}

\begin{Example}\label{ex:vnpi}
Suppose $\pi$ is a countable discrete group.  The group ring
$\C[\pi]$ acts on the left on $\ltwo(\pi)$.  The von Neumann
algebra $\VN{N}(\pi)$ is the von Neumann algebra generated
 by $\C[\pi]$, that is, the weak closure of $\C[\pi]$.
For $a = \sum_{g\in\pi}a_{g}\cdot g \in \VN{N}(\pi)$, the trace is
given by
\[ \Tr_{\pi}(a) = a_{1_{\pi}}. \]
\end{Example}

\subsection{Hilbert modules}
\begin{Def}
A \define{Hilbert module} for $\A$ is a Hilbert space $M$ with a
continuous left $\A$-module structure, and such that there exists a
Hilbert space $H$ and
an $\A$-equivariant isometric embedding $\iota$
of $M$ onto a closed subspace of the Hilbert space
$\hat\A\ctensor H$.
Here $\ctensor$ denotes the completed tensor product,
and $\A$ acts trivially on $H$.
Also note that $H$ and $\iota$ must
exist, but are not part of the structure.
\end{Def}

$B^{\A}(\hat\A\ctensor H)$ is identified with the algebra
$\A^{op}\ctensor B(H)$.
  If $(x_{1},x_{2},\dotsc)$ is an orthonormal
basis for $H$, $T \in B^{\A}(\hat\A\ctensor H)$ is identified with
a matrix with entries in $\A$ acting by right multiplication.
More precisely, let $P_{j} : \hat\A\ctensor H \to \hat\A$ be the
$j^{\text{th}}$ coordinate projection, given by
$P_{j}(a\tensor y) = \inner{y}{x_{j}}a$.  Then
\[
   T_{ij}(a) = P_{i}T(a\tensor x_{j}) : \hat\A \to \hat\A.
\]

Say that $T$ is \define{$\A$-Hilbert-Schmidt} if
$\sum_{i,j}\Tr_{\A}(T_{ij}^{*}T_{ij}) < \infty$. Say
that $T$ is of \define{$\A$-trace class} if
$T$ is the product of two $\A$-Hilbert-Schmidt operators.
If $T$ is $\A$-trace class,
the von Neumann trace of $T$ is
$\sum_{i=1,2,\dotsc}\Tr_{\A}T_{ii}$.
As usual, this does not depend on the basis of $H$.

If $M$ is a Hilbert module then orthogonal projection onto $M$,
$P_{M} : \hat\A\ctensor H \to \hat\A\ctensor H$, commutes with
$\A$.  Let
\[
  \dim_{\A}M = \Tr_{\A}P_{M}.
\]
This is well defined, \ie\ $\dim_{\A}M$ is independent of the choice
of $H$ and the embedding $\iota$.
For a bounded operator
$T : M \to M$ commuting with $\A$, define
$\Tr_{\A}T = \Tr_{\A}( P_{M}\circ T\circ\iota_{M})$
whenever the right hand side is of $\A$-trace class.

It is not hard to show that if $P : M \to M$ is a projection
and $T$ is bounded, then
\begin{equation}\label{eq:normtrace}
 \abs{\Tr_{\A}PT} \leq \norm{T}\dim_{\A}\Image(P),
\end{equation}
where $\norm{T}$ is the operator norm.
In particular if $\dim_{\A}M < \infty$ then for any bounded
$T : M \to M$, $\Tr_{\A}T$ is defined and in fact
\begin{equation}
 \abs{\Tr_{\A}T} \leq \norm{T}\dim_{\A}M.
\end{equation}

\begin{Example}
Let $\A = \C \subset B(\C,\C)$ as in Example~\ref{ex:triv}.
Then $\A$-Hilbert modules are simply $\C$-vector spaces,
we arrive at the standard definitions of trace class, and
$\Tr_{A}$ and $\dim_{\A}$ are the ordinary trace and dimension.
\end{Example}

% Just some explaination
\begin{comment}
the above is proved in your email bryan:
Suppose P :H->H  is self-adj, P^2 = P, and  Im P = M is finite dimensional.
Suppose T is any bounded operator.  Then  | Tr(PTP) | < ||T|| * dim M

Proof: (< always means \leq, all inner products and norms are H-S, except
||T|| at the end)

| Tr(PTP) | = | < TP, P > |  <  ||TP|| * ||P|| < ||T|| * ||P||^2 = ||T|| *
dim M

The first equality is the def. of the H-S inner product.
The next is Cauchy-Schwarz
The next is that ||TP|| < ||T||*||P|| for P H-S and T bounded.
Then for any projection,  ||P||^2 = Tr(P*P)  = Tr(P) = dim M.

\end{comment}
% end of explaination

%
\subsection{Determinants}
\label{sec:dets}
This section defines a notion of determinant for Hilbert modules.
Let $M$ be a finite dimensional Hilbert $\A$-module.
The determinant is defined for formal power series in $u$
with coefficients
in $\mathbf{B}_{M} = B^{\A}(M)$, using methods similar
to~\cite{Ba_zeta}.  These  determinants are formal power series over
$\C$. In our applications the series converge for small $u \in \C$ to
a complex-valued 
determinant which generalizes the real-valued Fuglede-Kadison
determinant defined in~\cite{fk:det}.

Let $u$ be an indeterminant.
For a ring $R$ let $R[[u]]$ denote 
the ring of formal power series in $u$ with coefficients in $R$,
and let $R^{+}[[u]]$ denote those power series with zero constant term.

Because $M$ is finite dimensional,
\[
 \Tr_{\A}: \mathbf{B}_M\to \field C
\]
is well defined and satisfies \eqref{eq:normtrace} for all
$T \in \mathbf{B}_{M}$.
The trace on $\mathbf B_{M}$ extends to
$\Tr_{\A} : \mathbf B_{M}[[u]] \to \C[[u]]$
by applying to each coefficient.

The exponential
$\Exp: \mathbf B_{M}[[u]]\to\Id+\mathbf B_{M}^+[[u]]$
and the principal branch of logarithm
$\Log:\Id +\mathbf B_{M}^+[[u]]\to \mathbf B_{M}[[u]]$
are defined as series.
Here $\Id$ is the identity operator on $M$.
$\Exp$ and $\Log$ are mutually inverse (see \cite{Ba_zeta} and
references therein).
Define
$\Det_{\A}: \Id+\mathbf B_{M}^+[[u]] \to \Id+\C^+[[u]]$
so that the following diagram is commutative:
\[
\begin{CD}
\Id+\mathbf B_{M}^+[[u]]      @>\Det_{\A}>>     1+\field{C}^+[[u]]\\
@V{\Log}VV                                  @V{\Log}VV\\
\mathbf B_{M}^+[[u]]          @>\Tr_{\A}>>     \field{C}^+[[u]]
\end{CD}
\]

% This is now partly subsumed by the general statment and partly
% included in the Prop:
%   Let $\mathbf B_{H_1,H_0}[[u]]$ be the set of power series in $u$ with
%   coefficient in $B^{\A}(H_0,H_1),$ and  $\mathbf B^+_{H_0,H_1}[[u]]$ be
%   the power series with $0$ for the constant term.

\newcommand{\boldal}{\boldsymbol{\alpha}}
\newcommand{\boldbe}{\boldsymbol{\beta}}
\begin{Prop}\label{det_pro}
The determinant $\Det_{\A}$ satisfies the following:
\begin{enumerate}
\item For $\mathbf{c} \in 1+\C^+[[u]]$,
      $\Det_{\A}(\mathbf{c}\Id) = \mathbf{c}^{\dim_{\A}M}
      = \Exp(\dim_{\A}M\Log \mathbf{c})$.
\item For $\boldal \in \Id+\mathbf B_{M}^+[[u]]$
      and invertible  $\boldbe\in \mathbf B_{M}[[u]]$,
      $\Det_{\A}(\boldbe \boldal \boldbe^{-1}) = \Det_{\A}(\boldal)$.
\item For $\boldal, \boldbe
           \in \Id+\mathbf B_{M}^+[[u]]$,
      $\Det_{\A}(\boldal\boldbe) =
       \Det_{\A}(\boldal) \Det_{\A}(\boldbe)$.
\item Let $H_0, H_1$ be finite dimensional Hilbert $\A$-modules.
      Then $B^{\A}(H_0 \oplus H_1)$ can be
      identified with matrices in the form
\[
  \alpha=
  \begin{bmatrix}
    \alpha_{00} & \alpha_{01}\\
    \alpha_{10} & \alpha_{11}
  \end{bmatrix} \in
  \begin{bmatrix}
    B^{\A}(H_0,H_0) & B^{\A}(H_1,H_0)\\
    B^{\A}(H_0,H_1) & B^{\A}(H_1,H_1)
  \end{bmatrix}.
\]
Let $\boldal\in \Id+\mathbf B_{H_0\oplus H_1}^+[[u]]$ be
\[
  \boldal =
  \begin{bmatrix}
    \boldal_{00} &      *      \\
         0       & \boldal_{11}
  \end{bmatrix}
\]
% I took this out and changed \alpha_{01} to * because
% alpha01 is a power series with coefficients in something that's
% not a ring, so it shouldn't be written as B[[u]]. --Bry
%  \in
%  \begin{bmatrix}
%    \mathbf B_{H_0}[[u]] & \mathbf B_{H_1,H_0}[[u]]\\
%    0 & \mathbf B_{H_1}[[u]]
%  \end{bmatrix}
with $\boldal_{ii} \in \mathbf B_{H_i}[[u]]$.
Then
\[
  \Det_{\A}(\boldal) =
  \Det_{\A}(\boldal_{00})\Det_{\A}(\boldal_{11}).
\]
\end{enumerate}
\end{Prop}

\begin{proof}
\begin{enumerate}
\item $\Log(\Det_{\A}(\mathbf{c}\Id)) = \Tr_{\A}\Log(\mathbf{c}\Id) =
       \Tr_{\A}(\Id\Log(\mathbf{c})) = \dim_{\A}M\Log(\mathbf{c})$.
\item $\Log\Det_{\A}(\boldbe \boldal \boldbe^{-1}) =
       \Tr_{\A}(\boldbe \Log(\boldal) \boldbe^{-1}) =
       \Tr_{\A}\Log(\boldal) =
       \Log\Det_{\A}(\boldal)$.
\item This follows from the Campbell-Baker-Hausdorff formula and the
      fact that $\Tr_{\A}$ vanishes on commutators.
\item This follows from
      $\Tr_{\A}\alpha = \Tr_{\A}\alpha_{00}+\Tr_{\A}\alpha_{11}.$
\end{enumerate}
See \cite[I-4]{Ba_zeta} for details.
\end{proof}

\begin{Prop}\label{prop:detconv}
Fix $T \in B^{\A}(M)$. Suppose $u \in \C$ satisfies
$\abs{u} < 1/\norm{T}$.
\begin{enumerate}
\item The series $\Det_{\A}(I - uT)$ converges absolutely.
\item $\abs{\Det_{\A}(I - uT)} = \Det_{\A}^{FK}(I-uT)$, where the
      right hand side is the Fuglede-Kadison~\cite{fk:det} determinant.
\end{enumerate}
% The proof of part 2 is not completely justified, but cleaning it up
% would obscure the simplicity and also add like a page of detail. Bry.
\end{Prop}
\begin{proof}
\[
  \Det_{\A}(I - uT) =
    \Exp\left(-\sum_{n=1}^{\infty}\frac{1}{n}\Tr_{\A}(uT)^{n}\right).
\]
The sum converges uniformly
since by \eqref{eq:normtrace},
$\abs{\Tr_{\A}(uT)^{n}} \leq \norm{uT}^{n}\dim_{\A}M$,
and since $\norm{uT} < 1$.

Now, recall that if $X \in B^{\A}(M)$ is invertible, the
Fuglede-Kadison determinant of $X$ is a positive real defined by
\[
  \log \Det_{\A}^{FK}(X) = \frac{1}{2}\Tr_{\A}\log(X^{*}X),
\]
where $\log(X^{*}X)$ is defined using the functional calculus for
self-adjoint operators.
Now we have
\begin{equation*}\begin{split}
  \log\abs{\Det_{\A}(I-uT)}^{2}
  &= \log\Det_{\A}(I-uT)\overline{\Det_{\A}(I-uT)}\\
  &= \log\Det_{\A}(I-uT)\Det_{\A}(I-\bar{u}T^{*})\\
  &= \log\Det_{\A}(I-uT)(I-\bar{u}T^{*})\\
  &= \Tr_{\A}\log(I-uT)(I-\bar{u}T^{*})\\
  &= \log (\Det_{\A}^{FK}(I-uT))^{2}.
\end{split}\end{equation*}
\end{proof}

\subsection{Hilbert representations}
Let $\G$ be any countable discrete group (although one could
allow more general groups), and $\A$ a finite von Neumann algebra.

\begin{Def}
A \define{Hilbert $\A$-representation} of $\G$ consists of
a Hilbert $\A$-module $M$ and a homomorphism
$\rho : \G \to B^{\A}(M)$.
\end{Def}

Say that $\rho$ is \define{unitary} if
$\rho(g)$ is a unitary
operator for all $g \in \G$.  The direct sum
of two Hilbert representations is defined in the usual sense.
One can form the tensor product over $\C$ of a Hilbert representation
with a finite dimensional representation of $\G$ to obtain a new
Hilbert representation.

For $G < \G$, let
\[
  M^{G} = \{ x \in M \mid \rho(g)x = x, \forall g \in G\}
\]
denote the invariant subspace of $M$ under $G$.  Since for $a \in
\A$, $g \in G$, and $x \in M^{G}$ we have $\rho(g)ax = a\rho(g)x = ax$,
$M^{G}$ is a Hilbert submodule of $M$.

\subsection{Group actions}
\label{sec:hrho}
Suppose now that $\G$ acts on a discrete countable set $X$, and
suppose the stabilizers $\G_{x}$ are finite for all $x \in X$.
Set $B = \G\bs X$.
Let $\rho : \G \to B^{\A}(M)$ be a unitary Hilbert representation.
In this section, we define a Hilbert $\A$-module
analogous to the space of $M$-valued $\ltwo$ functions on $B$, but
which is more natural when $\G$ does not act freely.

For $\G$-equivariant $f : X \to M$ define
\[
  \normrho{f}^{2} = \sum_{x \in B}\norm{f(\cover{x})}_{M}^{2}
\]
where for each $x$ we choose any lift $\cover{x} \in X$.
The definition is independent of the choice of lifts because
$\rho(\gamma)$ is unitary for all $\gamma \in \G$.

Now define
\[
  \HrhoX
  = \{ f: X \to M \mid
               f\text{ is $\G$-equivariant, and }
               \normrho{f} < \infty \}.
\]
$\HrhoX$ is a Hilbert space with norm $\normrho{\cdot}$ and inner
product
\[  \innerrho{f}{g} =
      \sum_{x \in B} \inner{f(\cover{x})}{g(\cover{x})}_{M}. \]

% \begin{Remark}
% The restriction to unitary representations is not essential, but
% suffices for the examples considered in this paper.  One could allow
% any uniformly bounded (in norm) representation, or if $B = \G\bs X$ is
% finite then any representation.
% In these cases, the topology on $\HrhoE$ is
% well defined, but the norm will no longer be canonical.
% \end{Remark}

$\HrhoX$ admits an $\A$ action given by
\[
  (a \cdot f)(x) = a \cdot f(x)
\]
where $a \in \A$, $f \in \HrhoX$, and $x \in X$.  An easy check shows
that the action preserves $\G$-equivariance of $f$.

Now choose a fixed lift $\cover{x} \in X$ for each $x \in B$.
Define $\Phi : \HrhoX \to M\ctensor\ltwo(B)$ by
\[ \Phi(f) = \sum_{x\in B}f(\cover{x})\tensor\delta_{x} \]
where $\delta_{x} \in \ltwo(B)$ is 1 at $x$ and 0 elsewhere.
Since the terms $f(\cover{x})\tensor\delta_{x}$ are mutually
orthogonal, $\Phi$ is an isometric embedding.
This shows that $\HrhoX$ is a Hilbert $\A$-module.  One should
note that the embedding $\Phi$ depends on the choice of lifts.

\begin{Prop}
The von Neumann dimension of $\HrhoX$ is given by
\begin{equation}\label{eq:dimHrho}
  \dim_{\A}\HrhoX = \sum_{x\in B}\dim_{\A}M^{\Gx}.
\end{equation}
\end{Prop}
\begin{proof}
Fix $x \in B$, and
recall that $\Gx$ denotes the stabilizer in $\G$ of
$\cover{x} \in X$.
For $\phi \in M^{\Gx}$, there is a function $f_{\phi} \in \HrhoX$
supported only on the orbit of $\cover{x}$,
\begin{equation}\label{eq:phifphi}
  f_{\phi}(\gamma \cover{x}) = \rho(\gamma)\phi
\end{equation}
It is easy to check that the correspondence
$\phi \leftrightarrow f_{\phi}$ induces an isomorphism of $\A$-Hilbert
modules
\[
  \hat{\bigoplus_{x\in B}} M^{\Gx}
     \stackrel{\sim}{\rightarrow}
  \HrhoX.
\]
This proves the Proposition, because the trace on $\A$ is normal.
\end{proof}

%
% SECTION:
% zeta-sec
%
%
%
\section{Zeta Functions}
\label{sec:zeta}
\subsection{Tree lattices}
Let $X$ be a locally finite tree, $\G$ an $X$-lattice (or a discrete
subgroup of $\Aut(X)$ without inversion).
Put $B = \G\bs X.$ Notation  introduced in this section will be
used  throughout the rest of this paper.\par

We denote the path metric on $X$ by $d.$ It is clear that $\G$ induces
an action on the ends of $X$, where an \define{end} is a class of
reduced infinite paths in X which eventually coincide.

For $g\in\G$ let $\len(g)=\min_{x\in X}d(gx,x)$
and let $\axis{g}=\{x\in X\mid d(gx,x)=\len(g)\}.$
Notice that $\axis{g}$ is the smallest subgraph of $X$ which is
invariant under the action of $\left< g \right>.$ If $\len(g)=0$
($g$ is \define{elliptic}), then $g$
is of finite order and $\axis{g}$ is the tree of fixed points of $g.$
We call $g$ \define{hyperbolic} if
$\len(g)>0.$ In this case $\axis{g}$ is an infinite ray, isomorphic to
$\Z,$ which we call the \define{axis} of $g.$ The induced
action of $g$ on $\axis{g}$ is just a translation of amplitude $\len(g)$
towards one of the two  ends of $\axis{g}$,
denoted by $\ep(g).$ When $g$ is hyperbolic $\ep(g)$ and $\ep(g^{-1})$
are the only two ends of $X$ fixed by $g.$

Let
\begin{equation}\label{equ:def:ends}
  \Ends = \Ends(\G) = \{\ep(g) \mid g\in\G,\len(g)>0\}
\end{equation}
be the set of hyperbolic ends of $X$, and put
$\Prim = \Prim(\G) = \G\bs\Ends(\G)$.

For $\ep\in\Ends(\G)$ let $\Gep$ be the stabilizer of $\ep.$
Let $[x,\ep)$ be the infinite ray from $x$ towards the end $\ep.$
Define $\tau_\ep: \Gep \to \Z$ by
\begin{equation}\label{equ:def:tau}
  \tau_\ep(\g) = d(x,u) - d(\g x,u)
\end{equation}
where
$x\in X$ and $u\in [x,\ep)\cap [gx,\ep).$ The function $\tau_\ep$ is a
homomorphism~\cite{Ba_zeta}.
For $\g\in \Gep,$ $\abs{\tau_\ep(\g)}=\len(\g).$
In particular, if $\ep=\ep(g)$ then $\tau_\ep(g)=\len(g).$ Put
\[
  \Gepo = \Ker(\tau_\ep)
\]
and define $\lenep>0$ so that
\[
  \lenep\Z = \Image(\tau_\ep).
\]

For $\ep=\ep(g)\in\Ends(\G)$ let $\axis{\ep} = \axis{g}$.
This definition is independent of $g$ because
all hyperbolic elements associated to the
same end have the same axis~\cite[II-2.4]{Ba_zeta}.
The graph $\axis{\ep}$ is invariant under $\Gep$, $\Gepo$ is finite and
acts trivially on $\axis{\ep}$, and the graph
$\Gep\bs\axis{\ep}$ is a cycle of length $\lenep.$

\subsection{Zeta functions of uniform tree lattices}
\label{sec:zetaul}
This section describes the zeta function associated to a uniform
lattice, due to Bass~\cite{Ba_zeta}.

Suppose $\rho:\G\to GL(V_\rho)$ is a finite dimensional
$\C$-representation of $\G.$
For $\ep\in \Ends(\G)$ put
\begin{equation}\label{eq:sigepdef}
  \sigep = \frac{1}{\abs{\Gepo}}
           \sum_{\substack{g\in\Gep \\ \tau_\ep(g)=\lenep}} g
\end{equation}

The zeta function associated to the action of $\G$ on $X$ and the
representation $\rho$ is defined as
\begin{equation}
  Z_\rho(u) = \prod_{\ep\in \Prim}
              \Det\left(\rho(1_\G-\sigep u^{\lenep})\right).
\end{equation}

For this section only,
let $\HrhoE$ (respectively $\HrhoV$) be the space of
$\G$-equivariant functions
from $EX$ (respectively $VX$) to $V_\rho.$
Define the operators $T_\rho :\HrhoE \to \HrhoE$ and
$\delta_\rho,Q_\rho : \HrhoV \to \HrhoV$ by
\begin{equation}\label{eq:defTdQ}
\begin{align}
    T_\rho f(e)      &= \sum_{\text{$(e,e_1)$ reduced}} f(e_1)\\
    \delta_\rho f(x) &= \sum_{\dz e=x}f(\doz e)\\
    Q_\rho f(x)      &= q(x)f(x)
\end{align}
\end{equation}
where $q(x)+1$ is the degree of $x$.
The operator $\delta_{\rho}$ is known as the \define{adjacency
operator}.

Bass' theorem~\cite{Ba_zeta} says that $Z_\rho(u)=\Det(I-T_\rho u)$,
which results in a  generalization of Ihara's theorem:
\[
  Z_\rho(u) = (1-u^2)^{-\chi_\rho(\G\bs X)}
              \Det(I-\delta_\rho u+ Q_\rho u^2),
\]
where
\[
  \chi_\rho(\G\bs X)=\sum_{x\in\G\bs X}\dim V_\rho^{\Gx}
       -\frac{1}{2}\sum_{e\in \G\bs EX}\dim V_\rho^{\Ge}
\]
is the Euler characteristic of $\G\bs X$ corresponding to $\rho.$

\subsection{Zeta functions of general tree lattices}
\label{sec:def_zeta}
Suppose $X$ is a locally finite tree with
a uniform bound $d$ on the degree of vertices.
In this section we define the zeta function associated to an
$X$-lattice $\G$
and a Hilbert representation of $\G$
which satisfies certain finiteness conditions. In fact the theory
applies to any discrete group $\G<\Aut(X)$ without inversion, equipped
with a representation which satisfies the finiteness
conditions. However, we 
concentrate only on the case where $\G$ is a lattice.

Let $\A$ be a finite von Neumann algebra with
trace $\Tr_{\A} : \A \to \C$.
Let $M$ be any Hilbert $\A$-module, and
$\rho : \G \to B^{\A}(M)$ a unitary Hilbert representation.

Using Section~\ref{sec:hrho}, we have Hilbert $\A$-modules
$\HrhoE$ and $\HrhoV$.  By \eqref{eq:dimHrho},
\begin{equation*}\begin{align}
  \dim_{\A}\HrhoE &= \sum_{e\in EB}\dim_{\A}M^{\Ge}\\
  \dim_{\A}\HrhoV &= \sum_{x\in VB}\dim_{\A}M^{\Gx}.
\end{align}\end{equation*}
Using these spaces,
define the operators $T_\rho,Q_\rho$ and $\delta_\rho$
as in \eqref{eq:defTdQ}.

Define
\[
  \Ends_n = \{\ep\in\Ends \mid \len(\ep) \text{ divides }n\}.
\]
The action of $\G$ on $\Ends$ preserves length, so $\G$ also acts
on $\Ends_{n}$.  Set $\Prim_n = \Prim_n(\G)=\G\bs\Ends_n.$

For $\ep \in \Ends$,
let $\axep$ denote the set of edges of $\axis{\ep}$ which are
oriented towards $\ep$.  Then define
\begin{equation*}
\begin{align*}
  \Endsn{E} &= \{(\ep,e) \mid \ep \in \Ends_{n}, e \in \axep\}\\
  \Endsn{e} &= \{\ep \mid \ep \in \Ends_{n}, e \in \axep\}.
\end{align*}
\end{equation*}
There is a natural diagonal action of $\G$ on $\Endsn{E}$, where
$\gamma\in\G$ sends $(\ep,e)\in\Endsn{E}$ to $(\gamma \ep,\gamma e).$

\begin{Lemma}\label{lem:dimHpEn}
With $\HrhoNE$ defined as in Section~\ref{sec:hrho},
\[
  \dim_{\A}\HrhoNE = \sum_{\ep\in \Prim_n} \lenep \dim_{\A}M^{\Gepo}.
\]
\end{Lemma}
\begin{proof}
Stabilizers in $\G$ of $(\ep,e) \in \Endsn{E}$ must fix $\ep$ and must
have translation length $0$ to fix $e$.
That is, $\G_{(\ep,e)} = \Gepo$.  Using \eqref{eq:dimHrho},
\begin{equation*}\begin{split}
  \dim_{\A}\HrhoNE
    &= \sum_{x \in \G \bs \Endsn{E}}\dim_{\A}M^{\Gx}\\
    &= \sum_{\ep \in \Prim_n}\sum_{e \in \Gep\bs\axep}
         \dim_{\A}M^{\G_{(\ep,e)}}\\
    &= \sum_{\ep \in \Prim_n}\lenep \dim_{\A}M^{\Gepo}.
\end{split}\end{equation*}
\end{proof}

\begin{Def}\label{def:fc}
Define two finiteness conditions on $\G$, $X$, and $\rho$:
\begin{enumerate}
\item $\dim_{\A}\HrhoE < \infty$.
\item For all $n$, $\dim_{\A}\HrhoNE < \infty$.
\end{enumerate}
\end{Def}
Note that $\dim_{\A}\HrhoE < \infty$ implies
$\dim_{\A}\HrhoV < \infty$, because $\HrhoV$ can be identified
with the subspace of $\HrhoE$ of
functions which are constant on edges that share initial vertices.

\begin{Remark}
Finiteness condition 2 is equivalent to the condition:
\begin{enumerate}
\item[$2^{\prime}$.] For all $n$,
  \begin{equation}\label{eq:finloop}
    \sum_{\substack{\ep \in \Prim;\\ \lenep = n}}
        \dim_{\A}M^{\Gepo} < \infty
  \end{equation}
\end{enumerate}
\begin{proof}
$\Endsn{E}$ is naturally
a disjoint union over $k|n$ of subsets $\{(\ep,e) \mid \lenep = k\}$.
Then $\HrhoNE$ splits as a finite direct sum over $k|n$, and the
dimension of the $k^{th}$ summand is exactly $k$ times the left hand
side of \eqref{eq:finloop}.
\end{proof}
\end{Remark}

Both finiteness conditions hold when $\G$ is a
uniform tree lattice and $\dim_{\A}M < \infty$.  This follows
easily from the fact that $\G \bs EX$ and $\G \bs \Endsn{E}$
are finite sets.
If $\G$ is a non-uniform $X$-lattice, the conditions
will hold if $M$ and $\rho$ are chosen appropriately.  See
Section~\ref{sec:nonuni} for more details.

\begin{Def}
Suppose $X$, $\G$, and $\rho$ satisfy the two finiteness conditions
of Definition~\ref{def:fc}.
Define $\sigep$ as in \eqref{eq:sigepdef}, and
the zeta function associated to $X$, $\G$, and $\rho$ as
\begin{equation}\label{eq:zetadef}
  Z_\rho(u) = \prod_{\ep \in \Prim}
              \Det_{\A}(I-\rho(\sigep)u^{\lenep})\quad\in \C[[u]]
\end{equation}
\end{Def}

In contrast with the zeta function in~\cite{Ba_zeta}, there may very
well be infinitely many $\ep \in \Prim$ with the same length.  This
means one must check that the coefficient
of $u^n$ is finite for each $n$, in order that
$Z_\rho(u) \in \C[[u]]$.  This question is answered positively in
the proof of Theorem~\ref{bass_hash}.

\subsection{Proof of the Bass-Hashimoto formula}
In this section, we prove Theorem~\ref{thm:mainZT}.
Formally, the proof is similar to
the arguments of~\cite{Ba_zeta}.
However there are
a number of convergence questions which must be addressed.
\newcommand{\sumep}{\sum_{\ep \in \Prim(\G)}}

\begin{Thm}\label{bass_hash}
Suppose $X$, $\G$, and $\rho$ satisfy the two finiteness conditions:
\begin{enumerate}
\item $\dim_{\A}\HrhoE < \infty$.
\item For all $n$, $\dim_{\A}\HrhoNE < \infty$.
\end{enumerate}
Then $Z_{\rho}(u)$ is well defined and
\begin{equation}\label{eq:ihara}
   Z_{\rho}(u) = \Det_{\A}(I-T_{\rho}u)\quad \in\C[[u]].
\end{equation}
In addition, the power series
$Z_{\rho}(u)$ converges for $u \in \C$, $\abs{u} < \norm{T_{\rho}}$.
\end{Thm}

\begin{proof}
Convergence for small $u$ will follow from \eqref{eq:ihara} and
Proposition~\ref{prop:detconv}.  So, compute formally:
\begin{equation}\begin{split}\label{eq:lhsofi}
\Log \Det_{\A}(I-T_{\rho}u)
  &= \Tr_{\A} \Log (I-T_{\rho}u) \\
  &= -\sum_{n=1}^{\infty}\frac{u^{n}}{n}Tr_{\A}T_{\rho}^{n}\\
  &= -\sum_{n=1}^{\infty}\frac{u^{n}}{n}
                \left({\sum_{\ep \in \Prim_{n}(\G)}}\lenep
                \Tr_{\A}\rho(\sigep^{n/\lenep})\right),
\end{split}\end{equation}
from Lemma~\ref{lem:tracet}.  The coefficient of $u^{n}$ converges
absolutely for each $n$.

Rearranging terms and substituting $n = m\lenep$,
\eqref{eq:lhsofi} becomes
\begin{equation*}\begin{split}
\Log \Det_{\A}(I-T_{\rho}u)
  &= -\sumep \sum_{m=1}^{\infty} \frac{u^{m\lenep}}{m}
                \Tr_{\A}\rho(\sigep^{m})\\
  &= \sumep \Tr_{\A} \left(-\sum_{m=1}^{\infty}
                \rho(\sigep^{m})\frac{u^{m\lenep}}{m}\right)\\
  &= \sumep \Tr_{\A} \Log\left(I - \rho(\sigep)u^{\lenep}\right)\\
  &= \sumep \Log \Det_{\A}\left(I - \rho(\sigep)u^{\lenep}\right)\\
  &= \Log Z_{\rho}(u).
\end{split}\end{equation*}
\end{proof}

%Computation of trace of T^n %
\begin{Lemma}\label{lem:tracet}
For all $n = 1,2,\dotsc$,
\begin{equation}\label{eq:trace_T^n}
  \Tr_{\A}T_{\rho}^{n} = {\sum_{\ep \in \Prim_{n}(\G)}}
                         \lenep\Tr_{\A}\rho(\sigep^{n/\lenep}).
\end{equation}
The sum converges absolutely.
\end{Lemma}
\begin{proof}
First we prove that the right hand side of \eqref{eq:trace_T^n} is
absolutely convergent. Notice that $\rho(\sigep^{n/\lenep})$ is a
projection to $M^{\G^0_{\ep}}$ followed by a unitary operator, hence
by~\eqref{eq:normtrace}
\[
  \abs{\Tr_{\A}\rho(\sigep^{n/\lenep})} \le \dim_{\A}M^{\G^0_{\ep}}.
\]
Using Lemma~\ref{lem:dimHpEn}
\[
  {\sum_{\ep \in \Prim_{n}(\G)}}
  \lenep \abs{\Tr_{\A}\rho(\sigep^{n/\lenep})}
  \le \dim_{\A}\HrhoNE
\]
which is finite by assumption.

As in Section~\ref{sec:hrho},
fix a lift $\cover{e} \in EX$ of each edge $e \in EB$.  We then
have an $\A$-equivariant isometric embedding
$\HrhoE \hookrightarrow \hat{\bigoplus}_{e \in EB}M$.  To compute
the trace of $T_{\rho}^{n}$, we calculate the trace of $T_{\rho}^{n}$
acting on the $e^{\text{th}}$ summand for each $e \in EB$.
More precisely, for $e \in EB$ let $(T_{\rho}^{n})_{e}$ denote
the composition
\begin{equation}\label{eq:seq}
  M \rightarrow M^{\Ge} \hookrightarrow \HrhoE
    \xrightarrow{{T_{\rho}^{n}}}
    \HrhoE \rightarrow M.
\end{equation}
The first map in (\ref{eq:seq}) is the projection
given by averaging over $\Ge$.
The second map is $\psi \rightarrow f_{\psi}$ described
in \eqref{eq:phifphi}.  The final map to $M$ is evaluation
at $\cover{e}$.  Thus
\[
  (T_{\rho}^{n})_{e} =
    \left(T_{\rho}^{n}f_{\rho(\delta_{\cover{e}})\phi}\right)
    (\cover{e}),
\]
where
\[
  \delta_{\cover{e}} = \frac{1}{\abs{\Ge}}\sum_{g\in\Ge} g.
\]
Using the definition of $T_{\rho}$, for any $\phi \in M$ we have
\newcommand{\nredsum}{\sum_{(\cover{e},e_{2},
                             \ldots,e_{n})\text{ red.}}}
\begin{equation*}\begin{split} %\allowdisplaybreaks
  (T_{\rho}^{n})_{e}\phi
  &= \nredsum f_{\rho(\delta_{\cover{e}})\phi}(e_{n})\\
  &= \nredsum \begin{cases}
       f_{\rho(\delta_{\cover{e}})\phi}(\gamma\cover{e})
       & \text{if $e_{n} = \gamma\cover{e}$ for
               some $\gamma \in \G$,}\\
       0 & \text{otherwise.}
     \end{cases}\\
  &= \nredsum \begin{cases}
       \rho(\gamma)\rho(\delta_{\cover{e}})\phi
       & \text{if $e_{n} = \gamma\cover{e}$ for
               some $\gamma \in \G$,}\\
       0 & \text{otherwise.}
     \end{cases}\\
  &= \nredsum \begin{cases}
       \rho(\gamma\delta_{\cover{e}})\phi
       & \text{if $e_{n} = \gamma\cover{e}$ for
               some $\gamma \in \G$,}\\
       0 & \text{otherwise.}
     \end{cases}\\
  &= \nredsum \rho\left(\frac{1}{\abs{\Ge}}
       \sum_{\gamma \in \G, \gamma\cover{e} = e_{n}}
       \gamma\right)\phi\\
  &= \sum_{\substack{\gamma \text{ hyp.}; \len(\gamma) = n;\\
                     \cover{e} \in X^{+}_{\gamma}}}
       \frac{1}{\abs{\Ge}}\rho(\gamma)\phi.
\end{split}\end{equation*}
We can now take the trace, to obtain
\begin{equation}\label{eq:tratn1}
  \Tr_{\A}(T_{\rho}^{n})_{e} =
     \sum_{\substack{\gamma \text{ hyp.}; \len(\gamma) = n;\\
                     \cover{e} \in X^{+}_{\gamma}}}
       \frac{1}{\abs{\Ge}}\Tr_{\A}\rho(\gamma).
\end{equation}

Bass~\cite[II.2.4(5)]{Ba_zeta} shows that all
hyperbolic elements associated to an end have the same axis.
We can therefore change (\ref{eq:tratn1}) to a sum over ends.

We have
\begin{equation}\begin{split}\label{eq:tratn2}
\Tr_{\A}(T_{\rho}^{n})_{e}
  &= \sum_{\substack{\ep \in \Ends_{n};\\ \cover{e}\in\axep}}
     \sum_{\substack{\gamma \in \Gep;\\ \tau(\gamma) = n}}
     \frac{1}{\abs{\Ge}}\Tr_{\A}\rho(\gamma) \\
  &= \sum_{\substack{\ep \in \Ends_{n};\\ \cover{e}\in\axep}}
     \frac{\abs{\Gepo}}{\abs{\Ge}}
     \Tr_{\A}\rho\left( \frac{1}{\abs{\Gepo}}
        \sum_{\substack{\gamma \in \Gep;\\ \tau(\gamma) = n}}
        \gamma \right) \\
  &= \sum_{\substack{\ep \in \Ends_{n};\\ \cover{e}\in\axep}}
     \left[ \Ge : \Gepo \right]^{-1}
     \Tr_{\A}\rho\left( \sigep^{n/\lenep} \right).
\end{split}\end{equation}
The trace of $T_{\rho}^{n}$ is the sum over $EB$ of
$\Tr_{\A}(T_{\rho}^{n})_{e}$.  That is, we have expressed the
trace as a sum over orbit representatives of edges and then a
sum over all ends.  The key step of this computation is to
switch to a sum over orbit representatives of ends and then a sum
over all edges.

Note that each $\Ge$ orbit of $\ep \in \Endsn{\cover{e}}$ has
$\left[ \Ge : \Gepo \right]$ elements.  Also, the final summand
in (\ref{eq:tratn2}) is $\G$-invariant as a function of $\ep$.
Thus,
\begin{equation}\begin{split}\label{eq:switch}
\Tr_{\A}T_{\rho}^{n}
  &= \sum_{e \in EB}
     \sum_{\ep\in\Endsn{\cover e}}
     \left[ \Ge : \Gepo \right]^{-1}
     \Tr_{\A}\rho( \sigep^{n/\lenep} )\\
  &= \sum_{e \in EB}
     \sum_{\ep \in \Ge \bs \Endsn{\cover{e}}}
     \Tr_{\A}\rho( \sigep^{n/\lenep} )\\
  &= \sum_{(e,\ep) \in \G\bs\Endsn{E}}
     \Tr_{\A}\rho( \sigep^{n/\lenep} )\\
  &= \sum_{\ep \in \G\bs\Ends_{n}}
     \sum_{e \in \Gep\bs\axis{\ep}}
     \Tr_{\A}\rho( \sigep^{n/\lenep} )\\
  &= \sum_{\ep \in \Prim_{n}} \lenep
     \Tr_{\A}\rho( \sigep^{n/\lenep} ).
\end{split}\end{equation}
The interchange of summations is allowed because the last sum
converges absolutely.
\end{proof}

\subsection{The operator $\Delta_\rho(u)$}
We prove Theorem~\ref{thm:mainLap} from the introduction.
Continue with the notation of  Section~\ref{sec:def_zeta}, and
suppose that $\dim_{\A}\HrhoE<\infty$
(and therefore $\dim_{\A}\HrhoV<\infty$).
\begin{Thm}\label{thm:main2}
Put $\chi_\rho(X)=\dim_{\A}\HrhoV - \half\dim_{\A}\HrhoE$.  Then
\begin{equation}
  \Det_{\A}(I-T_\rho u) =
  (1-u^2)^{-\chi_\rho(X)}\Det_{\A}(I -  \delta_\rho u+  Q_\rho u^2)
\end{equation}
where exponentiation is defined using the principal branch of the
logarithm.
\end{Thm}
\begin{proof}
The proof is very similar to that of~\cite[Theorem II.1.5]{Ba_zeta}.

Define the operators:
\begin{equation*}
\begin{align*}
        \dz,\doz &:\HrhoV \to \HrhoE \\
        \partial_i f (e)& =f(\partial_i e);\quad i=0,1 \\
        \sigma_0 &: \HrhoE \to \HrhoV \\
        \sigma_0 g (x)&= \sum_{\dz e=x}g(e)\\
        J&: \HrhoE\to\HrhoE\\
        Jg(e)&=g(\bar e)
\end{align*}
\end{equation*}
Notice that $\sigma_0\dz=I_V+Q_\rho$ and $\sigma_0\doz=\delta_\rho$.

Let $H=\HrhoV\bigoplus \HrhoE$ and define the following two operators on
$H$:
\[
L=\begin{bmatrix}
  (1-u^2) I_V  &  0 \\
  u\dz - \doz  &  I_E
  \end{bmatrix}, \qquad
M=\begin{bmatrix}
  I_V           &  u \sigma_0 \\
  u\doz - \dz   & (1-u^2)I_E
  \end{bmatrix}
\]
A simple computation shows that
\begin{equation*}
ML=
\begin{bmatrix}
\Delta_\rho(u) &  \sigma_0 u\\
       0       & (1-u^2)I_E
\end{bmatrix}
\end{equation*} and
\begin{equation*}
LM=
\begin{bmatrix}
(1-u^2) I_V & (1-u^2) \sigma_0 u \\
0 &           (I_E - T_\rho u)(I_E -Ju)
\end{bmatrix}.
\end{equation*}

$M$ is invertible (\cite[II-1.6]{Ba_zeta}) so
$\Det_{\A}(LM) = \Det_{\A}(M(LM)M^{-1}) = \Det_{\A}(ML)$
by Proposition~\ref{det_pro}.2.

Using Proposition~\ref{det_pro}, compute
\begin{equation*}
\begin{align*}
  \Det_{A}(ML)
    &= \Det_{\A}(\Delta_\rho (u))\Det_{\A}(1-u^2)I_E \\
    &= (1-u^2)^{\dim_{\A}\HrhoE}\Det_{\A}(\Delta_\rho (u))
\intertext{and}
  \Det_{\A}(LM)
    &= \Det_{\A}(I-T_\rho u)(I_E-Ju)\Det_{\A}((1-u^2)I_V) \\
    &= \Det_{\A}(I-T_\rho u)\Det_{\A}(I_E-Ju)\Det_{\A}((1-u^2)I_V) \\
    &= (1-u^2)^{\dim_{\A}\HrhoV}\Det_{\A}(I-T_\rho u)\Det_{\A}(I_E-Ju) \\
    &= (1-u^2)^{\dim_{\A}\HrhoV+\dim_{\A}\HrhoE/2}\Det_{\A}(I-T_\rho u).
\end{align*}
\end{equation*}
The last equality is the computation
$\Det_{\A}(I_E-Ju) = \half\dim_{\A}\HrhoE$ in~\cite[II-1.6]{Ba_zeta}.
The theorem follows easily.

\end{proof}

%
% SECTION:
% exam
%

%
%
\section{Coset Representations}
\label{sec:repex}
The coset representation of an infinite group on the Hilbert space of
$\Ltwo$ functions on a quotient group is the motivating example for
this paper.  In this section, the theorems of Section~\ref{sec:zeta}
are interpreted for coset representations.

Suppose there is an exact sequence $1\to \Lambda \to \G \to \pi \to 1$.
Take $\A = \VN{N}(\pi)$ (see Example~\ref{ex:vnpi}),
$M = \ltwo(\pi)$, and
$\rho$ the coset representation of $\G$ on $M$:
\[ (\rho(\gamma)\phi)(x) = \phi(\gamma^{-1} x) \]
for $\gamma \in \G$, $x \in \pi,\phi \in \ltwo(\pi).$
This situation will hold for the entire section.

\subsection{The spaces $H_{\rho}$}
Suppose $\G$ acts on a countable discrete set $X$, $B = \G\bs X$,
and set $Y = \Lambda\bs X$.
We have the natural action of $\pi$ on $Y$ and
$\pi\bs Y=B.$
For $y \in Y$, let
$\pi_{y} = \{\alpha \in \pi \mid \alpha y = y\}$.

\begin{Prop}\label{prop:dimpi}
The von Neumann dimension of $\HrhoX$ is given by
\[
  \dim_{\VN{N}(\pi)}\HrhoX =
  \sum_{x \in B}\frac{1}{\abs{\pi_{\hat{x}}}}
\]
where $\hat{x} \in Y$ is any lift of $x$.
\end{Prop}
\begin{proof}
Fix $x \in B$, with lift $\cover{x} \in X$ and $\hat{x} \in Y$.
From \eqref{eq:dimHrho}, we need to compute
$\dim_{\pi}\ltwo(\pi)^{\Gx}$ for each $x \in B$.
First notice that $\Gx/(\Gx \cap \Lambda) \cong \pi_{\hat{x}}$.

The action of $\G$ on $\ltwo(\pi)$ factors through $\pi$, and so
\begin{equation*}
\begin{split}
\ltwo(\pi)^{\Gx} &= \ltwo(\pi)^{\pi_{\hat{x}}} \\
                 &= \{f \in \ltwo(\pi) \mid
                      f \text{ is constant on cosets of }
                            \pi_{\hat{x}}\}.
\end{split}
\end{equation*}
Orthogonal projection $P$ onto $\ltwo(\pi)^{\Gx}$ is therefore
given by averaging over $\pi_{\hat{x}}$, and so
\( \inner{P(\delta_{1_{\pi}})}{\delta_{1_{\pi}}} =
   \abs{\pi_{\hat{x}}}^{-1} \).
\end{proof}

\begin{Remark}\label{rem:hrl2}
Since $\Lambda$ acts trivially on $\ltwo(\pi)$, functions in
$\HrhoX$ are periodic in the $\Lambda$ direction, so one can
think of them as functions on the quotient $Y$.
In fact, when $\Lambda$ acts freely on $X$
it is not hard to check that $\HrhoX$ and $\ltwo(Y)$
are isometric as $\VN{N}(\pi)$ modules by the map
$\Upsilon : \HrhoX \to \ltwo(Y)$, where
\[
  (\Upsilon f)(x) = \abs{\pi_{\hat{x}}}^{\half}f(\cover{x})(1_\pi).
\]
\end{Remark}

\subsection{The zeta function}
Now suppose that $X$ is a tree and $\G$ is an $X$-lattice.
Put $Y = \Lambda \bs X$.
As before, we let $\Ends(\G)$ denote the set of
hyperbolic ends of $\G$, and $\Prim(\G) = \G\bs\Ends(\G)$.
We also denote the hyperbolic ends of $\Lambda$ by $\Ends(\Lambda)$
and set $\Prim(\Lambda) = \Lambda\bs\Ends(\Lambda)$.
There is a natural action of $\pi$ on $\Prim(\Lambda)$  and we have
$\pi\bs\Prim(\Lambda) = \G\bs\Ends(\Lambda)$.

As before, we denote the length of an end $\ep \in \Ends(\G)$ by
$\len(\ep)$.
Since we have an injection $\Ends(\Lambda) \hookrightarrow \Ends(\G)$,
some ends of $\G$ are also ends of $\Lambda$ and therefore have a
possibly greater length.  To emphasize the difference, we
put
\[ \len_{\Lambda}(\ep) =
      \min\{\len(g) \mid
          g\ep = \ep; g \text{ hyperbolic}; g \in \Lambda\}
\]

\begin{Thm}\label{thm:zcomp}
When $Z_{\rho}(u)$ is well defined,
\begin{equation}\label{eq:zcoset}
  Z_{\rho}(u) = \prod_{\ep \in\pi\bs\Prim(\Lambda)}
                 \left(1-u^{\len_{\Lambda}(\ep)}\right)
                 ^{\frac{1}{\abs{\pi_{\ep}}}}
\end{equation}
where
$\pi_{\ep} = \{\alpha \in \pi \mid
             \alpha \cdot \ep = \ep \in \Prim(\Lambda) \}$.
\end{Thm}
\begin{proof}
From \eqref{eq:zetadef} we have
\[
  Z_\rho(u) =
    \prod_{\ep \in \Prim(\G)}
      \Det_{\pi}(I-\rho(\sigma_\ep)u^{\lenep})
\]
so that
\begin{equation}\label{eq:zrone}\begin{split}
  \Log Z_\rho(u) &=
     \sum_{\ep \in \Prim(\G)}
       \Tr_{\pi}\Log (I-\rho(\sigma_\ep)u^{\lenep})\\
  &= -\sum_{\ep \in \Prim(\G)}\sum_{n=1}^{\infty}
       \frac{u^{n\lenep}}{n}\Tr_{\pi}\rho(\sigma_\ep^{n}).
\end{split}\end{equation}

Now fix $\ep \in \Prim(\G)$.  Fix $n \in \N$.
Choose any $\gamma \in \Gep$ with $\len(\gamma) = n\lenep$.
We have
\begin{equation}\begin{split}
  \Tr_{\pi}\rho(\sigma_\ep^{n})
  &= \Tr_{\pi}\rho\left(\frac{1}{\abs{\Gepo}}\sum_{g \in \Gepo}
                        \gamma g\right)\\
  &= \frac{1}{\abs{\Gepo}}\sum_{g \in \Gepo}
     \begin{cases}
       1&\text{if $\gamma g \in \Lambda;$}\\
       0&\text{if $\gamma g \notin \Lambda.$}
     \end{cases}
\end{split}\end{equation}
If $\ep \in \Prim(\G) - \pi\bs\Prim(\Lambda)$, no hyperbolic
element of $\Gep$ lies in $\Lambda$.  So, for such $\ep$,
$\Tr_{\pi}\rho(\sigma_\ep^{n}) = 0$.

Define $\Lambda_{\ep}^{0} = \Gepo \cap \Lambda$.
Then for $\ep \in \pi\bs\Prim(\Lambda)$,
\begin{equation}\label{eq:trsign}
  \Tr_{\pi}\rho(\sigma_\ep^{n}) =
  \begin{cases}
     \frac{\abs{\Lambda_{\ep}^{0}}}{\abs{\Gepo}}
       & \text{if $\len_{\Lambda}(\ep)$ divides $n\lenep;$}\\
     0 & \text{otherwise.}
  \end{cases}
\end{equation}

Now, putting $k = n\lenep/\len_{\Lambda}(\ep)$,
(\ref{eq:zrone}) becomes
\begin{equation}\label{eq:zrtwo}\begin{split}
  \Log Z_\rho(u)
  &= -\sum_{\ep \in \pi\bs\Prim(\Lambda)}\sum_{k=1}^{\infty}
       \frac{\lenep u^{k\len_{\Lambda}(\ep)}}{k\len_{\Lambda}(\ep)}
       \left(\frac{\abs{\Lambda_{\ep}^{0}}}{\abs{\Gepo}}\right)\\
 % \text{(using Lemma~\ref{lem:diagram})}
  &= -\sum_{\ep \in \pi\bs\Prim(\Lambda)}\frac{1}{\abs{\pi_{\ep}}}
       \sum_{k=1}^{\infty}\frac{u^{k\len_{\Lambda}(\ep)}}{k}\quad
       \text{(using Lemma~\ref{lem:diagram})}\\
  &= \sum_{\ep \in \pi\bs\Prim(\Lambda)}\frac{1}{\abs{\pi_{\ep}}}
       \Log\left(1-u^{\len_{\Lambda}(\ep)}\right)\\
  &= \Log\prod_{\ep \in \pi\bs\Prim(\Lambda)}
       \left(1-u^{\len_{\Lambda}(\ep)}\right)^{\frac{1}{\abs{\pi_{\ep}}}}
\end{split}\end{equation}
which proves the theorem.
\end{proof}

\begin{Lemma}\label{lem:diagram}
With the notation of Theorem~\ref{thm:zcomp},
\[
  \abs{\pi_{\ep}} = \frac{\abs{\len_{\Lambda}(\ep)} \abs{\Gepo}}
                         {\abs{\lenep} \abs{\Lambda_{\ep}^{0}}}.
\]
\end{Lemma}
\begin{proof}
Define $\pi_{\ep}^{0} = \Gepo/\Lambda_{\ep}^{0}.$  Then we have
a commutative diagram:
\[
\begin{CD}
      @.   1    @.   1    @.   1    @.    \\
 @.      @VVV      @VVV      @VVV      @. \\
 1 @>>> \Lambda_{\ep}^{0} @>>> \Gepo @>>> \pi_{\ep}^{0} @>>> 1\\
 @.      @VVV      @VVV      @VVV      @. \\
 1 @>>> \Lambda_{\ep}     @>>>  \Gep @>>> \pi_{\ep} @>>>  1\\
 @.      @VVV      @VVV      @VVV      @. \\
 1 @>>> \len_{\Lambda}(\ep)\Z @>>> \lenep\Z @>>>
         \Z_{\len_{\Lambda}(\ep)/\lenep}  @>>>  1 \\
 @.      @VVV      @VVV      @VVV      @. \\
      @.   1    @.   1    @.   1    @.    \\
\end{CD}
\]
It is not difficult to check that all rows and columns are exact,
and the lemma follows.
\end{proof}

\subsection{Infinite graphs}
Suppose that $\Lambda$ acts freely on $X$.  Then the finiteness
conditions (Definition~\ref{def:fc}) needed to define the zeta
function and prove Theorem~\ref{thm:mainZT} are always satisfied.
This is the content of the following theorem.

\begin{Thm}\label{thm:freewell}
Suppose $\G$ is an $X$-lattice where $X$ is a tree of bounded degree. Let
$1\to \Lambda \to \G \to \pi \to 1$, and let $\rho$ be the coset
representation of $\G$ on $\ltwo(\pi)$.
If $\Lambda$ is free, then
\begin{enumerate}
\item $\dim_{\VN{N}(\pi)}\HrhoE < \infty$.
\item For all $n$, $\dim_{\VN{N}(\pi)}\HrhoNE < \infty$.
\end{enumerate}
\end{Thm}
\begin{proof}
\begin{enumerate}
\item
For any $e \in EB$, let $\cover{e}$ be a lift of $e$ to $EX$, and
let $\hat{e}$ be the edge in $EY$ covered by $\cover{e}$.
We have $\Ge = \Ge/(\Ge \cap \Lambda) \cong \pi_{\hat{e}}$ because $\Lambda$
is free and $\Ge$ is torsion.
Then by Proposition~\ref{prop:dimpi},
$\dim_{\VN{N}(\pi)}\HrhoE$ is the edge-volume of
$\G\bs X$ (Definition~\ref{def:lattice}),
which is finite.
\item
By combining Lemma~\ref{lem:dimHpEn} with the computation in
Proposition~\ref{prop:dimpi},
\begin{equation*}
\begin{align}
  \dim_{\VN{N}(\pi)}\HrhoNE
    &= \sum_{\ep\in\Prim_n} \lenep \dim_{\VN{N}(\pi)}
           \ltwo(\pi)^{\Gepo} \nonumber \\
    &= \sum_{\ep\in\Prim_n} \lenep
           \frac{1}{\abs{\Gepo}} \nonumber \\
    &\le n \sum_{\ep\in\Prim_n}
           \frac{1}{\abs{\Gepo}}. \label{tofinite}
\end{align}
\end{equation*}
Let $d$ be the bound on the degree of vertices of $X$.
We claim that $\abs{\Gepo} \ge d^{-n}\abs{\Ge}$,
where $\cover{e} \in X_{\ep}$.
To see the claim, let $\ep = \ep(g)$.  Then the orbit of
$g\cover{e}$ under $\Ge$
has at most $d^{\lenep} \le d^{n}$ elements. The group $\Gepo$ is
the stabilizer in $\Ge$ of $g\cover{e}$, and this proves the claim.

Choose a lift $\cover{e}$ for each edge $e \in EB$. Each class
$\ep \in \Prim_{n}$ has a representative which
contains at least one $\cover{e}$.
On the other hand, each
edge $\cover{e} \in EX$ is contained in only finitely many
$\ep \in \Prim_n$ (a loose bound is $d^n$).
Then continuing \eqref{tofinite},
\[
 n \sum_{\ep\in\Prim_n} \lenep \frac{1}{\abs{\Gepo}}
  \le n d^{2n} \sum_{e\in EB}\frac{1}{\abs{\Ge}}
\]
which is finite.
\end{enumerate}
\end{proof}

Combining Theorem~\ref{bass_hash}, Theorem~\ref{thm:main2},
and Theorem~\ref{thm:zcomp},
one arrives at the formula
\begin{equation}\label{eq:coform}
  Z_{\rho}(u)
  =
  \prod_{\ep \in\pi\bs\Prim(\Lambda)}
                 \left(1-u^{\len_{\Lambda}(\ep)}\right)
                 ^{\frac{1}{\abs{\pi_{\ep}}}}
  =
  (1-u^2)^{-\chi_\rho(X)}\Det_{\A}(I -  \delta_\rho u+  Q_\rho u^2).
\end{equation}
The remainder of this section consists of interpretation of
\eqref{eq:coform}.

When $\Lambda$ is free, $\Prim(\Lambda)$ has a nice
geometric interpretation as primitive, tail-less, backtrack-less
loops in $Y$ (or equivalently free homotopy classes).
The length of the loop associated to $\ep$ is $\len_{\Lambda}(\ep)$.
$\pi$ acts on these loops by translation, so the product in
(\ref{eq:coform}) is over translation classes of loops.
The group $\pi_{\ep}$ is the subgroup of $\pi$ which fixes the loop
associated to $\ep$.

Since the preceding discussion reformulates the definition
of $Z_{\rho}(u)$ completely in terms of $Y$ and $\pi$, it makes sense
to refer to the \define{$\pi$-zeta function of the
infinite graph $Y$},
\begin{equation}
   Z_{\pi}(Y,u) = Z_{\rho}(u)
                = \prod_{\gamma \in \pi\bs P}
                  \left(1-u^{\len(\gamma)}\right)
                  ^{\frac{1}{\abs{\pi_{\gamma}}}}.
\end{equation}
Here $P = \Prim(\Lambda)$.

It remains to interpret
the $\Laplace(u)$ side of formula \eqref{eq:coform}.
There are isomorphisms of $\HrhoE$ with $\ltwo(EY)$ and
of $\HrhoV$ with $\ltwo(VY)$ as described in Remark~\ref{rem:hrl2}.
The operator $\delta_{\rho}$ becomes the usual adjacency operator on $Y$ and
$Q_{\rho}$ is multiplication at $x \in VY$ by $deg(x) - 1$.
Then $\Laplace(u)$ is just a weighted version of the classical
Laplace operator on $\ltwo\ 0$-chains of $Y$.

Finally, the Euler characteristic $\chi_{\rho}(X)$ is the $\Ltwo$
Euler characteristic $\chi^{(2)}(Y)$, as in~\cite{cg:l2coho},
and is easily computable.
If $\pi$ acts freely on $Y$, $\chi^{(2)}(Y) = \chi(B)$.
Otherwise, using Proposition~\ref{prop:dimpi},
\[
  \chi^{(2)}(Y) = \sum_{x \in VB}\frac{1}{\abs{\pi_{x}}} -
                  \half\sum_{e \in EB}\frac{1}{\abs{\pi_{e}}}.
\]
This expression
is also studied in~\cite{cg:l2coho} where it is referred to as
the virtual Euler characteristic of $B$.

Under the identifications in this section, \eqref{eq:coform}
becomes
\[
   Z_{\pi}(Y,u) = (1-u^{2})^{-\chi^{(2)}(Y)}\Det_{\pi}(\Laplace(u)).
\]
This completes the proof of Theorem~\ref{thm:mainInf}.

\section{Examples}
This section is concerned with the situation of Section~\ref{sec:repex}.
That is,
there is an exact sequence $1\to \Lambda \to \G \to \pi \to 1$, and
$Y =  \Lambda\bs X$.
We take $\rho$ to be the coset representation of $\G$ on
$\ltwo(\pi)$.

\subsection{The infinite grid}
We examine Theorem~\ref{thm:mainInf} for a concrete example.
Let $\pi = \Z \times \Z = \left<a\right> \times \left<b\right>$.
Let $Y$ be the Cayley graph of $\pi$ with generators $(a,0)$ and
$(0,b)$.  $Y$ is an infinite grid.
The formula gives a method for computing the number $N(L)$
of translation classes of tail-less, backtrack-less
primitive loops in $Y$ of length $L$.  These are exactly the
number of possible shapes of closed paths on a grid, as should be made
clear by the picture to follow.

In this example, $\pi$ acts freely on all loops (this is true
whenever $\pi$ is torsion-free), so
\[
   Z_{\pi}(Y,u) = \prod_{\gamma \in \pi \bs P}(1-u^{\len(\gamma)})
                = \prod_{l = 1}^{\infty}(1-u^{l})^{N(l)}.
\]
Notice that $N(l) = 0$ for $l$ odd ($Y$ is bipartite).  Then
\begin{equation}\begin{split}
  -\Log Z_{\pi}(Y,u)
    &= -\sum_{l = 1}^{\infty}N(l)\Log(1-u^{l})\\
    &= -\sum_{m=1}^{\infty}\sum_{l|m}N(l)\frac{l}{m}u^{m}\\
    &= -\sum_{M=1}^{\infty}\sum_{L|M}N(2L)\frac{L}{M}u^{2M}.
\end{split}\end{equation}
Computing coefficients of this power series will then compute
(inductively) the numbers $N(2L)$.  To achieve this, examine the
$\Laplace(u)$ side of the formula in Theorem~\ref{thm:mainInf}.  We have
\begin{equation}\begin{split}
  -\Log Z_{\pi}(Y,u)
    &= \Tr_{\pi}(-\Log(\Laplace(u)))
       + \sum_{M = 1}^{\infty}\frac{1}{M}u^{2M}\\
    &= \sum_{n=1}^{\infty}\Tr_{\pi}(1-\Laplace(u))^{n}
       + \sum_{M = 1}^{\infty}\frac{1}{M}u^{2M}\\
\end{split}\end{equation}

The weighted Laplacian $\Laplace(u)$ acts on
$\ltwo(VY) \cong \ltwo(\pi)$. We write it as a
$1 \times 1$ matrix with coefficient in $\Z[\pi]$, acting on
$\ltwo(\pi)$ by right multiplication:
\[
  \Laplace(u) = (1 - (a + a^{-1} + b + b^{-1})u + 3u^{2}).
\]
Computing $\Tr_{\pi}(1-\Laplace(u))^{n}$ is an exercise in
combinatorics.  One finds that
\[
  -\Log Z_{\pi}(Y,u) = \sum_{M = 1}^{\infty}
     \left[ \frac{1}{M}+\sum_{d = 0}^{M}\frac{(-3)^{M-d}}{M+d}
            \binom{M+d}{M-d}{\binom{2d}{d}}^{2} \right] u^{2M}.
\]

The first few values of $N(2L)$ are $0,2,4,26,152,1004,\dotsc$.
The 26 types of loops of
length 8 are pictured below (each has two possible orientations).

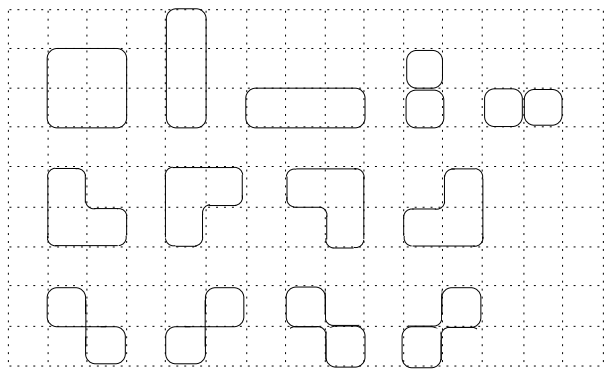
\picta{cm99_loop8grid}{8cm}{Loops of Length 8}

\subsection{Non-uniform tree lattices}
\label{sec:nonuni}
If $\G$ is a uniform $X$-lattice then any
Hilbert representation of $\G$ on a space of finite
von Neumann dimension satisfies the
finiteness conditions needed to define the corresponding zeta function.
In particular, one can let $\A = \C$ and choose any representation
of $\G$ on a finite dimensional vector space.
However, except for the trivial representation there is no
interesting canonical choice of representation.

When $\G$ is a non-uniform lattice the
trivial representation of $\G$
does not satisfy the finiteness condition---$\HrhoE$ is not
finite dimensional.
However in this case if
one can find a free normal subgroup $\Lambda \normal \G$
then the corresponding coset representation yields
a well-defined zeta function, by Theorem~\ref{thm:freewell}.

There are other possible representations for non-uniform lattices.
The following example is a coset representation with a well-defined
zeta function, where $\Lambda$ is not free.

\begin{Example}
Suppose $\G$ acts on a $q+1$-regular tree $X$
with quotient as in the diagram:

\picta{cm99_Padicex}{10cm}{The quotient $\G\bs X$}

To understand the figure, notice that
it shows only  geometric edges. For every oriented edge there is
an assigned index which is displayed  in the
figure at the initial vertex of the
edge.

For edge $e = (x_{n},x_{n-1})$
the index is $q$. Then if $\cover x_{n}$ is
a lift of $x_{n}$ there are exactly $q$ edges
$e_{1},\dotsc,e_{q}$ covering $e$ with
$\dz e_{i} = \cover x_{n}$.
The same interpretation holds for other indices.

Any such $\G$ is always an $X$-lattice (see \cite{bass_lattice}).

Choose a lift $(\cover x_{0}, \cover x_{1},\dotsc)$ of
the ray $\G\bs X$ to $X$.
Then $\G = \left< \G_{\cover x_n}; n=0,1,\dotsc \right>$, the group
generated by the vertex stabilizers.
In fact, $\G_{\cover x_n} < \G_{\cover x_{n+1}}$ for all $ n\ge 1$.

Let
$\Lambda =
     \left< \G_{\cover x_0},\G_{\cover x_1}\right>_{\text{normal}}$,
the normal subgroup generated by $\G_{\cover x_0}$ and
$\G_{\cover x_1}$ in $\G$.  Form the exact sequence
\[
  1\to \Lambda \to \G \to \pi \to 1
\]
where $\pi=\left< \G_{\cover n_i}; n=2,3,\dotsc \right>$.
It is not hard to show that $\dim_{\pi} \HrhoE < \infty$ if we assume
that $\G_1\normal\G_i,\ i\ge 1.$
We claim that $\Prim_n$ is finite for
each $n$ and therefore finiteness condition 2
(Definition~\ref{def:fc}) is satisfied.
The axis of each hyperbolic element of $\G$ must
pass through some lift of $x_0$. To count hyperbolic elements up to
conjugacy we need only count those whose axes pass through a fixed
lift of $x_{0}$, and there are only finitely many such axes.

We have shown that the zeta
function for $\G$ and the coset representation $\rho$
of $\G$ on $\Lambda$ is well defined.

Using the trivial representation for $\G$, the above
argument shows that finiteness condition 2 may be satisfied
although $\dim_{\A}\HrhoE$ is not finite.
\end{Example}

%
%  Bibliography
%

\end{document}